\documentstyle[12pt]{article}
\begin{document}
\newcommand{\be}{\begin{equation}}
\newcommand{\e}{\end{equation}}
\newcommand{\ba}{\begin{eqnarray}}
\newcommand{\ea}{\end{eqnarray}}
\newcommand{\n}{\nabla}
\newcommand{\dl}{\Delta}   
\newcommand{\x}{(x,t)}
\newcommand{\s}{S^n\times (0,T)}
\newcommand{\p}[1]{\frac{\partial #1}{\partial t}}
\newcommand{\k}[1]{\nabla_{#1}}
\newcommand{\kl}[1]{#1_{k l}}
\newcommand{\pl}[1] {\frac{\partial}{\partial t} #1}
\renewcommand{\theequation}{\arabic{section}.\arabic{equation}}
\newcommand{\no}{\noindent}
\newcommand{\nnb}{\nonumber}
\newcommand{\q}{equation}
\newcommand{\Q}{equations}

\title{Deforming Convex Hypersurfaces to a Hypersurface
with Prescribed Harmonic Mean
Curvature \thanks{ Supported by National Science 
Fundation of China and by grant of Morningside Center of Mathematics
,CAS}}
\author{{Huai-Yu Jian}\\
{\small Department of Applied Mathematics, Tsinghua University,
Beijing 100084, P.R.China}\\
{\small (e-mail: hjian@math.tsinghua.edu.cn)}}

\date {}
\maketitle

{\bf Abstract.} Let $F$ be a smooth convex and positive function
 defined in $A=\{x\in R^{n+1}: R_1 <|x|<R_2 \}$ satisfying
 $F(x)\geq nR_2$ on the sphere $|x|=R_2$ and $F(x)\leq nR_1$
 on the sphere $|x|=R_1$. In this paper, a heat flow method
 is used to deform convex hypersurfaces in $A$ to a hypersurface
 whose harmonic mean curvature is the given function $F$.

\vskip 0.4cm



\vskip 0.4cm

\section*{1. Introduction}
\setcounter{section}{1}
Let $M$ be a smooth embedded hypersurface in $R^{n+1}$ and $k_1,k_2,\cdots,k_n$ be its 
principal curvatures. Then $H^{-1}$ is called the harmonic mean curvature 
of $M$ if
\be H={\frac {1}{k_1}}+{\frac{1} {k_2}} +\cdots + {\frac{1} {k_n}}.\e
\no The question which we are concerned with is that given a function $f$
defined in $R^{n+1}$, under what conditions does the equation
$$H^{-1}(X)=f(X), X\in M$$
\no has a solution for a smooth,closed,convex and embedded hypersurface $M$, where
$X$ is a position vector on $M$.

The kind of such question was proposed by S.T.Yau in his famous problem
section [Y]. Many authors have studied the cases of mean curvature and Gauss
curvature instead. See, for instance, [BK], [TW], [T], [CNS], and [TS1] for the 
mean curvature and [O1], [O2], [TS2] for the Gauss curvature,  [CNS], [G1]
 and [G2]
for general curvature functions.

Let $F=f^{-1}$, then the problem above is equivalent to looking for a smooth,closed,
convex and embedded hypersurface $M$ in $R^{n+1}$ such that
\be H(X)=F(X), X\in M  \e
\no where $H$ is the inverse of harmonic mean curvature given by (1.1).
We are interested only in the hypersurfaces $M\subset A,$  a ring domain 
defined by
   $$A=\{ X\in R^{n+1} : R_1 < |X|<R_2\}.$$ 
for some constants $R_2 > R_1 >0.$   For this purpose, we need to suppose that 
$F$ is a smooth 
positive function defined in $R^{n+1}$ satisfying

\parindent 28pt {\bf(a)} $F(X)>nR_2$ for $|X|=R_2$ and $F(X)<nR_1$ for $|X|=R_1$,

\noindent and 

\parindent 28pt {\bf(b)} F is concave in A.

We will use a heat flow method to deform convex 
hypersufaces to a solution to (1.2). That is we conside the parabolic \q
\be \left.
\begin{array}{l}
\p X =\left (H(X)-F(X)\right )\nu (X), X\in M_t, t\in (0, T) \\
X(\cdot , 0) \ \ given, 
\end{array}
\right\} \e

\no where $X(x,t) $ : $S^n\rightarrow R^{n+1}$ is the parametrization of $M_t$
given by inverse Gauss map, which will be solved, and $\nu (X)$ is the outer
normal at $X\in M_t$, so $\nu (X\x )=x$ by the definition. Of course, $M_0$ 
is a given initial hypersurface.

\vskip 0.4cm
The following is our main result of this paper.

{\bf Theorem 1.1.} Suppose that $F$ is a smooth positive function satisfying
conditions (a) and (b), and a initial hypersurface $M_0 \subset A$ is smooth,
uniformly convex and embedded, satisfying $H(X_0 )\geq F(X_0)$ for all $X_0
\in M_0 .$ Then equation (1.3)
has a unique smooth solution for  $T=\infty $ which parametrizes a family
of smooth,closed,uniformly convex and embedded hpersufaces ,
$\{M_t : t\in [0, \infty )\}.$ Moreover, there exists
a subsequence $t_k \rightarrow
\infty $  such that  $M_{t_k}$ converges to a
smooth,closed,uniformly convex and
embedded hypersurface which lies in $\bar A $ and solves problem (1.2).

\vskip 0.4cm
We will have to meet two difficulties in proving the result above.
One is the gradient estimate for the support function, $u\x ,$ of the
hpersufaces $ M_t $
 which solves the equation (2.4) below; the other is the proof of convexity preserving
 (see the inequality (2.5)). In order overcome the first difficulty, we use 
 the well-known equality $|X\x |^2 = u^2 +|\n u|^2 $ and a general geometric 
 result suggested by R. Bartnik to estimate the $|X\x |^2$ for the equation (1.3).
 The proof of convexity-preserving is inspired by the computations in Hamilton
 [H1,H2] and Huisken [HU].

{\bf Remark 1.2.} Applying a usual approximation method, one can replace 
condition (a) in theorem 1.1 by

\parindent 28pt {\bf (a')} $ F(X)\geq nR_2$ for $|X|=R_2$ and $F(X)\leq R_1$
for $|X|=R_1.$

\vskip 0.3cm

After completing the paper, the author was aware that similar results had
been obtained by Gerhardt in [G1]. Although the results in [G1] are very
general, they don't include the above theorem 1.1 because our curvature
function $H$ does not belong the class $(K)$ in [G1] and our initial
hypersurfaces $M_0$ may be arbitrary instead of the fixed barrier
$M_1 =\{|X|=R_1\} $ in [G1]. Moreover, the arguements are absolutely different.

\vskip 0.3cm

{\small \it {\bf Acknowledgement.} The author would like to thank professors
R. Bartnik, K. Tso and W.Y. Ding for many helpful conservations and professor
B. Chow for his interests in this work.}

 \section*{2. Evolution equations, convexity-preserving and global solutions}
\setcounter{section}{2}
\setcounter{equation}{0}
 In this section, we will at first reduce the \q \ \ (1.3) to a equivalent 
 quasilinear parabolic \q \ \ on $S^n$ for the support function of $M_t ,$
 then we will show that this parabolic problem is globally solvable and it
 preserves convexity.

 We recall some facts in [U;p.97-98]. Let $e_1, e_2, \cdots , e_n$ be a 
 smooth local orthonormal frame on $S^n$, and let $\n _i=\n _{e_i}$,
 i=1, 2, ..., n and $\n =(\n _1, \n _2, \cdots , \n _n )$ be the covariant
 derivatives and the gradient on $S^n$, respctively. Since $X\x $ is the inverse
 Gauss map, the support function of $M_t$ is given by
 \be u\x =<x,X\x >, x\in S^n \e
 \no where $ <\cdot , \cdot >$ denotes the usual inner product in $R^{n+1}.$
 The second fundamental form of $M_t$ is
 $$h_{i j}(X\x )=\n _i \n _j u\x +\delta_{i j} u\x , i,j=1,2,\cdots , n.$$
 \no If $M_t$ is uniformly convex, then $h_{i j}$ is invertible, and hence
 the inverse harmonic mean curvature is the sum of all the eigenvalues of the matix
 $$b_{i j}=\left[ g^{i k}h_{k j}\right]^{-1}, $$

 \no where $g_{i j}$ is the metric of $M_t$. But the Gauss-Weingarten relation
 $$\n _i x =h_{i k}g^{k l}\n _l X $$
 and the fact $<\n _i x , \n _j x >=\delta_{i j}$ imply $g_{i j}
 =h_{i k}h_{j k}.$ Therefore, $b_{i j}=h_{i j}$ and 
 \be  H={\frac{1}{k_1}} +\cdots +{\frac{1}{k_n}}=\dl u +nu , \e

 \no where $\dl =\sum_{i=1}^{n}\n _i \n _i .$ Furthermore, since $x, \n _1 x , 
 \n _2 x, \cdots , \n _n x$ form a standard orthonormal basis at point $X\x ,$
 so we have
 \ba X\x & = & <X , x >x +<X , \n _i x>\n _i x \nnb  \\
 & = & u x + \n _i <X , x>\n_i x \nnb  \\
 & =& ux +\n _i u \n _i x .    \ea
 Using the results above and repeating the argument [U: p.98-101] in verbatim, 
 one has obtained the following lemma.

{\bf Lemma 2.1.} If for $t\in [0 , T)$ with $T\le \infty \ \ X\x $
 is a solution of (1.3) which parametrizes a smooth,closed,
 uniformly convex and embedded hypersurface , then the support
 functions $u\x $ of $M_t$ satisfy

 \be  \left.
 \begin{array}{l}
 \p u = \dl u +nu-F(ux+\n _i u\n _ix ), \ \ \x \in \s  \\
 u(\cdot , 0)=u_0 (\cdot ), \ \ the \ \ support \ \ function \ \ of \ \ M_0
 \end{array}
 \right\} \e

 \no and
 \be \n ^2 u +uI >0  \ \ in \ \ \s, \e
 \no where $I$ denotes the $n\times n $ unit matrix.  Conversely, if $u$ is a 
 smooth solution to (2.4) and satisfies (2.5), then the hypersurface $M_t$, determined by
 its support function $u\x $, is a smooth,closed,uniformly convex ,
 embedded hypersurface and solves (1.3) for
  $t\in [0, T).$

 From now on, we assume that the initial hypersurface $M_0$ is smooth,closed,uniformly
 convex . That is $u_0 \in
 C^\infty (S^n)$ and for some positive constant $C_0$,
 \be C_0 I\le \n ^2 u_0 +u_0 I . \e
 \no Noting that (2.4) is a quasilinear parabolic equation on the compact
 manifold $S^n$, by standard result for short-time existence (see, for example,
 [H3]), we have

{\bf Lemma 2.2.} There exists a maximal existence time $T=T(u_0)\in (0,\infty ]$
 such that (2.4) has a unique smooth solution  $u\in C^{\infty}(\s )
 \bigcap C([0,T); C^{\infty}(S^n)).$ If $T< \infty , $
 then 
 $$\lim_{t\rightarrow T^-}\|u(\cdot , t)\|_{C^1 (S^n)}= \infty . $$

{\bf Remark 2.3.} One can use the contraction principle and repeat the same argument 
 as in proving the short-time existence for harmonic heat flows to give a 
 direct proof of this lemma. See [ES] or [D].

 For the sake of deriving an apriori estimates, we need the following geometric
 result which was suggested by R. Bartnik.

{\bf Lemma 2.4.}  Let $X$ be the positive vector of a
smooth,closed hypersurface $M$ in $R^{n+1}$
 with outer normal $\nu (X)$ at $X\in M .$ Then if $|X|=<X , X>^{\frac{1}{2}}$
 attains a maximun $R$ at a point $X_0 \in M ,$ then $X_0 =R\nu (X_0)$
 and
 $$\amalg (w,w)\geq {\frac{1}{ R}} g(w,w), \forall w\in T_{X_0}M;$$
 \no if $|X|$ attains a minimum $r$ at a point $X_0 \in M ,$ then 
 $X_0 =r\nu (X_0)$ and
 $$\amalg (w,w)\leq {\frac{1}{ r}} g(w,w), \forall w\in T_{X_0}M_0 , $$
 \no where $g$ is the metric on $M$ and $\amalg$ is the second fundamental
 form  of $M$ with respect to the direction $-\nu .$

{\bf  Proof.} We consider only the first case, because the latter is completely
 anagolous. For $X_0\in M$ and any $w\in T_{X_0}M, $ choose a curve 
 $\gamma (s)$ on $M, \gamma : [0,1]\rightarrow M ,$ such that 
 $$\gamma (0) =X_0 , \ \ \ \  \dot {\gamma}=w . $$
 \no Let $\rho (X)=|X|. $ since $\rho ^2 (X)$ attains a maximum at 
 $X_0 \in M,$ then at this point
 $$\n \rho ^2 =0  \ \ and  \ \ \n ^2 \rho ^2 \leq 0.$$
 Therefore, we have
 $${\frac{d}{ds}}\rho ^2 (\gamma (0))=\n \rho ^2 (X_0)\cdot \dot {\gamma}(0) = 0$$
 and
 $${\frac{d^2}{ds^2}}\rho ^2 (\gamma(0))=\dot \gamma (0)\cdot \n ^2 \rho ^2 
 (X_0)\cdot \dot {\gamma}(0)+\n \rho ^2 (X_0)\cdot \ddot{\gamma}(0)
 \leq 0 . $$
 On the other hand,
 $${\frac {d \rho ^2}{ds}} = {\frac{d}{ds}}|\gamma (s)|^2 = 2\gamma (s)
 \cdot \dot \gamma (s)$$
 and
 $${\frac{1}{2}} {\frac{d^2 \rho ^2 }{ds^2}}=|\dot \gamma (s)|^2 +\gamma (s)
 \cdot \ddot \gamma (s). $$
 Thus,
 \be 0=2\gamma (0)\cdot \dot \gamma (0)= 2X_0\cdot w \e
 and
 \be 0\geq |\dot \gamma (0)|^2+X_0\cdot \ddot \gamma (0).  \e
 Since $w\in T_{X_0}M $ can be arbitrary, (2.7) implies
 $$X_0 =|X_0 |\nu (X_0)= R\nu (X_0), $$
 and so,
 \ba X_0\cdot \ddot \gamma (0) & = & R<\nu (X_0),\ddot \gamma (0)> \nnb  \\
 & = & R<\nu (X_0) , D_{\dot \gamma (0)} \dot \gamma (0)> \nnb  \\
 & = & -R\amalg (\dot \gamma (0), \dot \gamma (0)) \nnb  \\
& = &  -R\amalg (w,w), \nnb  \ea
which, together with (2.8), gives us that
\begin{eqnarray*}
\amalg (w,w) & \geq & {\frac{1}{R}} |\dot \gamma (0)|^2  \\
& = & {\frac{1}{ R}} <\dot \gamma (0), \dot \gamma (0)>   \\
& = & {\frac{1}{ R}} <D_{\dot \gamma (0)} X, D_{\dot \gamma (0)} X >  \\
& = & {\frac{1}{R}} g(\dot \gamma (0) , \dot \gamma (0))  \\
& = & {\frac{1}{ R}} g(w,w). 
\end{eqnarray*}

{\bf Lemma 2.5.} Suppose that $M_0 \subset \subset A $ in addition to (2.6). Let 
$u$ be a smooth solution to (2.4) and satisfies (2.5) on $\s $ with $T\leq 
\infty .$ Then for all $\x \in \s ,$ we have
$$ R_1^2< u^2 \x +|\n u\x |^2 < R_2^2. $$

{\bf Proof.} It follows from lemma 2.1 that the position vector $X(\cdot,t)$
of $M_t$ determined by the support function $u\x $ satisfies (1.3), i.e.
\be \left.
\begin{array}{l}
\p X =(H(X)-F(X))\nu (X), X\x \in M_t, \x \in \s  \\
X(\cdot , 0)=X_0 .
\end{array} \right\} \e
Moreover, (2.3) and the fact that $x, \n _1 x, \cdots , \n _n x$ form a 
standard orthonormal basis imply that
\be |X|^2=u^2+|\n u|^2 , \forall \x \in \s . \e
Thus, it is sufficient to prove that for $X$ solving (2.9),
\be R_1^2<|X\x |^2 <R_2^2, \forall \x \in \s . \e

For each $t\in [0, T),$ let
$$P_{min}(t)=\min_{X\in M_t}|X|^2 =\min_{x\in S^n }|X\x |^2, $$
and
$$P_{max}(t)=\max_{X\in M_t}|X|^2=\max_{x\in S^n}|X\x |^2.  $$
By virture of the assumption for $M_0,$ we have (2.6) and
\be R_1^2 <P_{min}(0)\leq P_{max}(0)<R_2^2. \e
Since $M_t$ is smooth, $P_{min}$ and $P_{max}$ are obviously Lipschitzian
on $[0,T).$  Were the inequality (2.11) not true, then by (2.12) we could find 
$t_1$ and $t_2$ in $[0,T)$ such that either
\be P_{min}(t_1)=R_1^2 \e
or
\be P_{max}(t_2)=R_2^2. \e
Without loss of generality, we assume that the case (2.13) happens, and the case
(2.14) is completely similar. Let 
$$t^* =inf \{t\in (0,T) : P_{min}(t)=R_1^2\},$$
and choose $x^*\in S^n $ such that
$$P_{min}(t^*)=|X(x^*,t^*)|^2 .$$
In order to compute the principal curvatures $k_1, \cdots , k_n$ of $M_{t^*}$
at $X(x^*,t^*),$ we use the principal direction $\xi _1 , \cdots , \xi _n $
to obtain
$$\amalg (\xi _i ,\xi _i )=k_i<\xi _i ,\xi _i>=k_i g(\xi_i,\xi_i),$$
so
$$k_i =\amalg (\xi _i,\xi _i)g^{-1}(\xi _i,\xi _i),   i=1,2,\cdots, n .$$
Thus, lemma 2.4 implies that $X(x^*,t^*)=R_1\nu (X(x^*,t^*))$ and
$$k_i\leq {\frac{1}{R}_1}, \ \ i=1,2,\cdots, n $$
at $X(x^*,t^*).$ Therefore, we have
$$H(X(x^*,t^*))\geq nR_1,$$
and 
\begin{eqnarray*}
\p {|X|^2} (x^*,t^*) & = & X\cdot \p X  \\
& = & R_1 \nu \cdot (H-F)\nu \nnb  \\
& \geq  & R_1 [nR_1 -F(X(x^*,t^*))]  \\
& >  0, 
\end{eqnarray*}
where we have used the condition (a) for $X(x^*,t^*)=R_1.$

On the other hand, we obviouly have
$$|X(x^*,t)|^2 > |X(x^*,t^*)|^2, \forall t \in [0, t^*) $$
thus
$$\p {|X|^2} (x^*,t^*)\leq 0$$
which yields a contradiction.

Next, we will prove that the convexity of $M_t$ is preserved. That is, (2.5) 
remains true for all $t\in (0,T)$ if it is so at $t=0.$

Let
$$\kl h =\k k \k l u +\kl {\delta} u, \ \ k, l=1,2,\cdots , n.$$
It follows directly from (2.4)and (2.3) that
\be \pl ({\kl {\delta}}u)=\kl {\delta}\left [ \sum_{i=1}^n h_{i i}-F(u x+X\x)\right 
] .\e
On the other hand, we differentiate the \q \ \ (2.4) twice to get
\be \pl {\k k \k l u}= \k k \k l {\dl u} +n\k k \k l u -F_j (X)\k k X^j
-F_{j h}(X)\k k X^h \k l X^j.   \e
Using the standard formula for interchanding the order of covaraint 
differentiation with respect to the othonormal frame on $S^n$, we have
\be \k k \k l \dl u =\dl \k k \k l u - 2n \k k \k l u +2 \kl {\delta}
\dl u,   \e
see, for instance, [CLT, p.85]. Thus, combining (2.15)-(2.17), we get
\ba  \pl {\kl h } & = & \dl \kl h +\kl {\delta} \dl u -n \k k \k l u
+\kl {\delta} \left [ \sum_{i=1}^n h_{ii}-F(X\x )\right ] \nnb  \\
& & - F_j (X\x ) \k k \k l X^j -F_{j h} (X)\k k X^h \k l X^j .  \ea
Set
$$ G\x = \sum_{i=1}^n h_{i i}\x -F(X\x ). $$
By virture of (2.18), (2.3), and (2.4), we see that
\ba \p G & = & \pl {\sum_{i=1}^n h_{ii}} - \p F   \nnb  \\
& = & \dl G +nG -\p F  \nnb \\
& = & \dl G +n G - F_j (X)(\p u x^j +\k i \p u \k i x^j ) \nnb  \\
& = & \dl G + n G -F_j (X) (Gx^j +\k i G \k i x^j ). \ea

{\bf Lemma 2.6.}  Suppose that in addition to (2.6), the support function $u_0$
of the initial hypersurface $M_0$ satisfies
\be G(x, 0)=\dl u_0(x) +n u_0(x) -F(u_0 +\k i u_0 \k i x )\geq 0, 
\ \ x\in S^n \e
and let $u\x $ be a smooth solution to (2.4) on $\s $ with $T\leq \infty .$
Then for all $\x \in \s , $
\be G\x = \dl u\x  + n u\x - F(ux +\k i u \k i x )\geq 0 . \e

{\bf Proof.} It is sufficient to prove that for each $T_0 <T,$ (2.21) holds true 
for all $\x \in S^n\times[0, T_0).$
Let
$$C_1 = max\left\{|u\x |+|\n u\x |: \x \in S^n\times [0, T_0]\right\}, $$
$$C_2 = max \left \{\sum_{j=1}^n |{\frac{\partial F}{\partial y_j}}|
:|y|\leq C_1\right\}, $$
and
$$G_{min}(t)=min\left \{ G\x : x\in S^n\right\}.$$
If $G_{min}(t)\leq 0 $ for some $t\in(0,T_0)$ we could find $x_t\in S^n$
such that
$$G(x_t ,t)=G_{min}(t)\leq 0, \ \ \dl G(x_t ,t)\geq 0, \ \ \n G(x_t ,t)=0 .$$
Thus, (2.19) implies that at $(x_t ,t)$
$$ \p G  \geq  n G -F_j (X)Gx^j \geq (n+C_2)G. $$
This yields
$$\pl {\left( Ge^{- (n+C_2)t}\right) }(x_t ,t)  \geq 0 .$$

Note that $e^{-(n+C_2)t}G_{min}(t)$ is a Lipchitizian function of $t .$
Denote it by $B(t)$ for simplicity. Then the result above gives us that
$$\liminf_{\tau \rightarrow 0^+}{\frac {B(t+\tau)-B(t)}{\tau}}\geq 0
\ \  when  \ \  B(t)\leq 0  \ \  for  \ \  t\in (0,T_0).$$
Also see [U, p.107]. Now by a result of Hamilton[H2, lemma 3.1] we 
conclude that
$$B(t)\geq 0, \ \ \forall t\in [0, T_0].$$
This proves the lemma.

{\bf Lemma 2.7.} Assume that the support function $u_0$ of the initial
hypersurface $M_0$
satisfies (2.6) and (2.20). Let $u\x $ be a smooth solution to (2.4) on
$\s $ with $T\leq \infty.$ Then for all $\x \in \s , $ we have
$$ \n ^2 u \x +u \x I >0. $$

{\bf Proof.} If the conclusion were not true, we could find a finite number
$t_0 \in (0, T)$ such that the minimum eigenvalue of the matrix
$[\kl h (x,t_0)]$
is zero, but  $[\kl h \x ]$ is positive definite for all $\x \in S^n\times
[0, t_0).$ Thus the inverse matrix $[h^{p q}\x ]$exists for all such $\x .$
Let
$$\kl g \x =<\k k X\x , \k l X\x >, \ \ k,l=1, 2, \cdots , n$$
be the metrix of $M_t.$ Then Gauss-weingarten relation and the fact
$\kl {\delta} =<\k i x , \k j x>$ gives us that
$$\k i x^j =h_{i k}g^{k l}\k l X^j, \ \ \ \ g_{i j}=h_{i k }h_{j k}.$$
See [U, p.98]. Therefore, for each $\x \in S^n\times [0, t_0),$ we have
\ba \k k \k l X^j & = & \k k \left( g_{l m}h^{m i}\k i x^j \right) \nnb  \\
& = & \k k \left(  h_{l i} \k i x^j \right) \nnb  \\
& = & \k k \k i x^j h_{i l}+\k i x^j \k k h_{l i}.
\nnb \ea
On the other hand, it easily follows from the standard formula for commuting
the order of covariant differentiation on $S^n$ that
$$\k i h_{j k}=\k j h_{i k}=\k k h_{i j} \ \ for \ \ all \ \ i, j=1, 2, \cdots ,
n $$
see [U, p.111]. With the aid of two equatities above, (2.18) turns to
\ba \pl {h_{k l}}& = & \dl \kl h +\kl {\delta}\dl u - n \k k \k l u +\kl{\delta}
G  \nnb  \\
& & -F_j (X)\left( \k k \k i x^j h_{i l}+\k i x^j \k i h_{l k} \right)
-F_{j h} (X)\k k X^h \k l X^j .  \ea
Now let us suppose that the minimum eigenvalue of $[\kl h ]$ over $S^n$
at time $t\in[0,t_0)$ attains at a point $x_t \in S^n.$ By rotating the 
frame $e_1, e_2, \cdots , e_n ,$ we may assume that $h_{1 1}(x_t,t)$ is
the minimum eigenvalue and $h_{j 1}(x_t,t)=0$ for $j=2, 3, \cdots ,n .$
(Also see [CLT, p.89].) Combining (2.22)and lemma 2.6 together, at 
$(x_t, t)$ we have
\be \p {h_{1 1}} \geq \dl u -n\k 1 \k 1 u -F_j(X)\k 1 \k 1 x^j h_{1 1}
-F_{j h}\k 1 X^j \k 1 X^h, \e
where we have used $\dl h_{1 1}\geq 0$ and $\nabla h_{1 1}=0$ at $(x_t , t).$
Since at $(x_t,t)$
$$nu +n\k 1 \k 1 u \leq \sum_{i=1}^n h_{i i}=\dl u +nu, $$
we have
$$\dl u -n\k 1 \k 1 u\geq 0. $$
Applying this equality and the concavity of $F$ to (2.23), we obtain
$$\p {h_{1 1}}\geq -F_j (X)\k 1 \k 1 x^j h_{1  1}.$$
Denote
$$C_3 =max\left\{ |u\x |+|\n u\x | :\x \in S^n\times [0, t_0]\right\},$$
$$C_4 =max\left \{ \sum_{j=1}^n |F_j (Y)| : |Y|\leq C_3 \right\}, $$ 
and 
$$C_5=C_4 max \left \{ |\k 1 \k 1 x | : x\in S^n \right \} . $$
Noting $h_{1 1}(x_t, t)$ is positive for all $t\in [0, t_0),$ we have
$$\pl {h_{1 1}(x_t ,t)} \geq -C_5 h_{1 1}(x_t ,t), \ \ t\in [0, t_0 ).$$
That is
$$\pl {\left ( h_{1 1}(x_t ,t)e^{C_5 t}\right )}\geq 0 , \ \  t\in [0, t_0).$$
By the maximum principle and the condition (2.6), we see that
$$h_{1 1}(x_t ,t)\geq h_{1 1}(x_0, 0)e^{-C_5 t}\geq C_0 e^{-C_5 t}  \ \
\forall t\in [0, t_0). $$
Letting $t\rightarrow t_0 ^- , $ we obtain
$$h_{1 1}(x_{t_0}, t_0)\geq C_0 e^{-C_5t_0 } , $$
which is contradictory with the assumption that the minimum eigenvalue of
$[\kl h ]$ is zero. This proves the lemma

{\bf Theorem 2.8.} Suppose that the initial hypersurface $M_0 \subset \subset A$ and its
support function $u_0$ satisfies (2.6) and (2.20). Then there exists a 
unique smooth solution $u$ to the following problem:
\be \left. 
\begin{array}{l}
\p u =  \dl u +n u -F(ux+\k i u \k i x) \ \ \  in \ \ S^n\times (0, \infty ) \\
u(\cdot , 0)=u_0(\cdot ) \ \ \ \  on \ \ S^n , \end{array}
\right\} \e
Moreover, for all $(x,t)\in S^n\times (0, \infty ), \ \ u\x $
satisfies
\be \n ^2 u +u I>0, \e
\be \p u =\dl u +nu -F(ux +\k i u \k i x)\geq 0, \e
and
\be R_1^2 <u^2 +|\n u|^2 <R_2^2 . \e
{\bf Proof.} By virture of lemma 2.2, we know that the problem (2.24) has a 
unique smooth solution on $\s $ with some $T\leq \infty . $ Moreover,
lemma 2.6 and 2.7 tell us that both (2.25) and (2.26) are satisfied on
$\s  .$ Thus, by lemma 2.5, we see that (2.27) is true for all $\x \in
\s .$ Using lemma 2.2 again, we know that $T$ is nothing but $\infty .$
This completes the theorem.

\section*{3. Converging to a convex hypersurface }
\setcounter{section}{3}
\setcounter{equation}{0}

In this section, we will use the theorem 2.8 in last section to prove 
the main result of this paper, theorem 1.1.

We begin with choosing a smooth,closed,uniformly convex hypersurface $ M_0 $
such that its support function
satisfies (2.6) and (2.20). (The existence of such $ M_0 $ is obvious
due to the condition (a)).
By theorem 2.8 and lemma 2.1
we obtain a family of smooth,closed,uniformly convex hpersufaces
$ M_t $ whose position vectors are
$$X\x =u\x x +\k i u\x \k i x , \ \ \x \in S^n\times [0, \infty), $$
where $u\x $ are the support functions of $M_t$ and satisfy
(2.24)-(2.27) in theorem 2.8. It follows from (2.26) and (2.24) that 
$$\left(\p u \right)^2 = \p u (\dl u +nu -F(X\x )\leq \p u (\dl u +nu).$$
Hence
$$\int_{S^n}(\p u)^2 dx \leq {\frac{1}{2}} \pl {\int_{S^n} (nu^2-|\n u|^2)dx} ,
\forall t\in (0, \infty ).$$
This, together with (2.27), implies that
\be \int_0^{\infty}\int_{S^n}(\p u )^2dxdt\leq (n+2)R_2^2\cdot vol (S^n). \e
Since (2.27) implies that for each $\alpha \in (0, 1)$
$$ nu -F(ux+\k i u\k i x)\in L^{\infty}(S^n\times [0, \infty )) ,$$
it follows from a property of heat \q \ \ that
\be ||u(\cdot ,t)||_{C^{1,\alpha }(S^n)}\leq C \e
uniformly in $t\in[1, \infty ).$ See [LSU, ch.4] or [D]. Moreover, by
(3.2)and the same argument, we have
\be ||u(\cdot ,t)||_{C^{3, \alpha}(S^n)}\leq C \e
uniformly in $t\in [1, \infty ). $ Using (3.1) and (3.3), we can find a 
sequence $t_k \rightarrow \infty$ as $k\rightarrow \infty $ such that
\be \int_{S^n}({\frac{\partial u}{\partial t_k}})^2dx \rightarrow 0 , \e
and
\be u(\cdot , t_k)\rightarrow U_0 \ \ in \ \ C^{3,\alpha}(S^n). \e
Furthermore, (2.24),(3.4) and (3.5) gives us
$$\int_{S^n}[\dl U_0 +nU_0 -F(U_0 x+\k i U_0 \k i x)]^2dx=0, $$
which yields
\be \dl U_0 +nU_0 -F(U_0 x+\k i U_0 \k i x)=0,
\ \ x\in S^n .\e
Since $U_0\in C^2(S^n),$ applying a standard elliptic theory to the \q
(3.6), we know $U_0 \in C^{\infty}(S^n).$

On the other hand, it easily follows from (3.4),(3.5) and (3.6)that as $k
\rightarrow \infty$
$${\frac{\partial u}{\partial t_k }}(\cdot , t_k)\rightarrow 0 \ \  in \ \
C(S^n).$$
Hence
\be   \n {\frac{\partial u}{\partial t_k}}(\cdot,t_k) \rightarrow 0 \ \
in \ \ C(S^n). \e
Therefore, (3.5)and (3.7) implies that
\be X(x,t_k)=u(x,t_k)x+\k i u(x,t_k)\k i x\rightarrow U_0 x+\k i U_0 \k i
x \ \ in \ \ C^1(S^n), \e
and
\be {\frac{\partial X}{\partial t_k}}(x,t_k)={\frac{\partial u}{\partial t_k}}
(x,t_k)x+\k i {\frac{\partial u}{\partial t_k}}(x,t_k)\k i x
\rightarrow 0 \ \ in \ \ C(S^n). \e
Since each hypersurface $M_t$ is smooth,closed,uniformly convex ,
lemma 2.1 says that $X\x $ satisfy (1.3)
for all $t\in [0, \infty ).$ Thus, using (2.2), (3.5), (3.8) and (3.9), we know that
\be Y_0(x)  =U_0 (x)x+\k i U_0 \k i x  \ \ \ \ \  x\in S^n  \e
satisfies the \q (1.2).

Now let $M_y$ be the hypersurface  determined by the position vector $Y_0(x), x\in
S^n.$ Obviously, $M_y$ is closed and encloses the origin, and is convex
because of (3.5), (3.8) and the uniformly convexity of hpersufaces $M_t$ for any $t\in
[0, \infty ).$ So $M_y$ is uniformly convex because $Y_0$
satisfies (1.2). Furthermore, due to the fact $U_0\in C^{\infty}(S^n), $
$M_y$ is smooth, too. Finally, using (2.27), (3.5) and (3.10), we obtain
$$R_1^2\leq |X_0(x)|^2\leq R_2^2, \ \ x\in S^n, $$
which means that the hypersurface $M_y$ lies in $\bar{A} .$ This prove theorem 1.1.

\newpage
\no {\bf References}

\no [BK] L.Bakelman and B. Kantor, Existence of spherically homeomorphic
 hypersurfaces in Euclidean space with prescribed mean curvature, Geometry
 and Topology, Leningrad(1974), 3-10.

 \no [CNS] L.Caffarelli, J. Nirenberg and J. Spruck, The Dirichlet problem
 for nonlinear second order elliptic \q . IV: Starshaped compact Weingarten
 hpersufaces , Current topics in partial differential \Q , ED. by Y. Ohya, K. Kasahara, N.
 Shimakura (1986), 1-26, Kinokunize Co., Tokyo.

 \no [CLT] B. Chow, L.P. Liou and D.H. Tsai, On the nonlinear prabolic \q \ \ 
 $\partial _t u = F(\dl u +nu)$ on $S^n$, Comm. Anal. Geom. 3(4)(1996), 75-94.

 \no [D] W.Y. Ding, Lecture Notes on harmonic maps, Inst. Math. Academia
 Sinica, 1995.

\no [G1]  C. Gerhardt, Closed weingarten hypersurfaces in Riemannian manifolds,
J. Differ. Geom. 43(1996), 612-641.

 \no [G2] C. Gerhardt, Hypersurfaces of prescribed weingarten curvature, 
 Math. Z. 224(1997), 167-194.

 \no [ES] J. Eells and J.H. Sampson, Harmonic mappings of Riemannian
 manifolds, Amer. J. Math. 86(1964), 109-160.

 \no [H1] R. S. Hamilton, Three-manifolds with positive Ricci curvature, J. Differ.
 Geom. 17(1982), 255-306.

 \no [H2] R. S. Hamilton, Four-manifolds with positive Ricci curvature, J. Differ.
 Geom. 24(1986), 153-179.

 \no [H3] R.S. Hamilton, Lecture Notes on heat \q \ \ in geometry, Honolulu,
 Hawaii, 1989.

 \no [HU] G. Huisken, Flow by mean curvature of convex hpersufaces into spheres,
 J. Differ. Geom. 20(1984), 237-266.

 \no [LSU] O. Ladyzeskaja, A. Solonnikov and N. N. Uralceva, Linear and
 quasilinear \Q \ \  of parabolic type, Trans. Amer. Math. Soc. 23, Providence, 1968.

 \no [O1] V.I. Oliver, Hypersurfaces in $R^{n+1}$ with prescribed 
 Gaussian curvature and related equations of Morge-Ampere type, 
 Comm. p.d.e., 9(1984), 807-839.

 \no [O2] V.I. Oliver, The problem of embeddeding $S^n$ into $R^{n+1}$
 with prescribed Gauss curvature and its solution by variational methods,
Trans. Amer. Math. Soc., 295(1986), 291-303.

 \no [T] A. Treibergs, Existence and convexity of hpersufaces of prescribed curvature,
 Ann. Scuola Norm. Sup.Pisa CI. Sci. 12(1985), 225-241.

 \no [TW] A. Treibergs and S.W. Wei, Embedded hpersufaces with prescribed mean curvature,
 J. Differ. Geom. 18(1983), 513-521.

 \no [TS1] K. Tso, On the existence of convex hpersufaces with prescribed mean curvature,
 Ann. Scuola Norm. Sup. Pisa 16 (1989), 225-243.

 \no [TS2] K. Tso, Existence of convex hpersufaces with prescribed Gauss-Kronecker
 curvature, J. Differ. Geom. (1991), 389-410.

 \no [U]J.I.E. Urbas, An expansion of convex hypersurfaces , J. Differ. Geom.
 33(1991), 91-125; Correction to, ibid. 35, 763-765.

 \no[Y] S.T. Yau, Problem section, Seminar on differential geometry, Ed. S.T. Yau,
 Ann. of Math. Stud., Princeton Univ. Press. 102, 1982, 669-706.

 \end{document}